




\documentstyle{amsppt}
\input amsppt1
\pageheight{19.8true cm}
\pagewidth{12.2true cm}

\parindent=4mm
\parskip=3pt plus1pt minus.5pt
\nologo\NoRunningHeads\NoBlackBoxes

\def\i{\looparrowright}

\def\f{\flushpar }
\def\nl{\newline }

\def\x{\times }

\def\f{\flushpar }
\def\nl{\newline }

\def\x{\times }

\document
\topmatter
\title
The projections of $n$-knots which are not \\
the projection of any unknotted knot    
\endtitle
\affil
Eiji Ogasa\\
ogasa\@ms.u-tokyo.ac.jp    
\endaffil
\thanks{
\f{\bf Mathematics Subject Classification (1991):} 57M25, 57Q45 
\nl This research was partially supported by Research Fellowships 
of the Promotion of Science for Young Scientists.  
\nl {\bf Keywords.} The projections of $n$-knots, Kummer surface, spun knots
}
\endtopmatter

\f{ \bf Abstract. }   
Let $n$ be any integer greater than two.   
We prove that there exists a projection $P$ 
having the following properties.
(1) $P$ is not the projection of any unknotted knot.   
(2) The singular point set of $P$ consists of double points.    
(3)$P$ is the projection of an $n$-knot which is diffeomorphic to the standard sphere.   

We prove 
there exists an immersed $n$-sphere ( $\subset\Bbb R^{n+1}\x\{0\}$ ) 
which is not the projection of any $n$-knot ($n>2$). 
Note that the second theorem is different from the first one. 

\baselineskip10pt
\head \S 1.   Introduction and Main results \endhead     
The study of   $n$-dimensional knots and links has a long history.
The research was connected with surgery theory.
(See 
[9], [20] 
etc. for the history. )
There are many fruitful results. There are many outstanding problems.
For example,  $n$-links have not been classified.  
 (This open problem is not discussed explicitly in this paper 
but it is one of motivations of this paper.) 

When we study $n$-knots and $n$-links ($n\geq2$), 
we sometimes see similarities and differences
between  the theory of 1-links and that of n-links. 
In this paper we point out one difference between them, 
associated with the projections of knots. 
 (Theorem 1.3.).

When one studies classical knots (in $\Bbb R^3$),   
it is important to consider the projections of knots into  $\Bbb R^2$.   
See  
[1], 
[3],  
[8], 
[12], 
[16], 
[23], 
[25], 
etc.   
For 2-dimensional knots in $\Bbb R^4$,   
one considers the projections of 2-knots into  $\Bbb R^3$.
See 
[2],  
[4], 
[5], 
[7],    
[10], 
etc. 

In order to state our problems (Problem 1.1 and 1.2) and our main theorem (Theorem 1.3.), 
we prepare some definitions on $n$-knots and on their projections.

We work in the smooth category.  

An {\it (oriented)  n-(dimensional) knot $K$} is 
a smooth oriented submanifold of $\Bbb R^{n+1}\x\Bbb R$ 
which is PL homeomorphic to the standard $n$-sphere.   
We say that n-knots $K_1$ and $K_2$ are {\it equivalent} 
if there exists an orientation preserving 
\nl diffeomorphism 
$f:\Bbb R^{n+1}\x\Bbb R\rightarrow\Bbb R^{n+1}\x\Bbb R$
\nl  such that $f(K_1)$=$K_2$  and 
$f\vert_{K_1}:$ $K_1$ $\rightarrow$ $K_2$ 
is an orientation preserving diffeomorphism.    
We say that n-knots $K_1$ and $K_2$ are {\it identical} 
if $id(K_1)$=$K_2$ holds for the identity map 
$id:\Bbb R^{n+1}$ $\x$ $\Bbb R$ $\rightarrow$ 
$\Bbb R^{n+1}\x\Bbb R$
and $id\vert_{K_1}:$ $K_1$ $\rightarrow$ $K_2$ is 
an orientation preserving diffeomorphism.    

Let $T$ be the unit $n$-sphere in 
$\Bbb R^{n+1}$ $\x$ $\{0\}$ $\subset$ $\Bbb R^{n+1}$ $\x$ $\Bbb R$.  
An $n$-knot $K$ is said to be {\it unknotted} 
if $K$ is equivalent to $T$.  

Let 
$\pi$:  $\Bbb R^{n+1}$ $\x$ $\Bbb R$ $\to$ $\Bbb R^{n+1}$ 
be the natural projection map.    
We assume $\pi\vert_K$ is a self-transverse immersion.  
The {\it projection} $P$ of an $n$-knot $K$ is 
$\pi\vert_K(K)$ in $\Bbb R^{n+1}$.    
We give $P$ an orientation by using 
the orientation of $K$.
The {\it singular point set} of 
the projection of an $n$-knot $K$ is 
the set  $\{$ $x\in$ $\pi\vert_{K}$ $(K)$  $\vert$  
$\sharp\{(\pi\vert_{K})$ $^{-1}(x)\}\geqq2.\}$.     

It is well-known that the projection of any 1-dimensional knot 
is the projection of an unknotted 1-knot.   
This fact is used in the definition of the Jones polynomial. 
It is also used in some definitions of the Conway-Alexander polynomial.      
See, e.g.,  
[13], 
and 
[14]. 

It is natural to ask the following question.   

\f{\bf Problem 1.1.} 
Let $P$ be the projection of an $n$-knot $K$. 
Let $K$ be diffeomorphic to the standard sphere.     
Is $P$ necessarily the projection of an unknotted $n$-knot? 

In this paper we give a negative answer to Problem 1.1 when $n>2$. 

The problem in the case $n=2$ is still open.   
As mentioned above, the answer in the case $n=1$ is affirmative.  

Note that, when $n=1$,  the singular point set always consists of double points.  

Here, we consider the following problem 1.2. 
Obviously, 
if the answer to Problem 1.2 is negative, 
the answer to Problem 1.1 is negative.  
Of course, as mentioned above, the answer to Problem 1.2 for $n=1$ is affirmative.

\definition {Problem 1.2} 
Let $P$ be the projection of an $n$-knot $K$. 
Let $K$ be diffeomorphic to the standard sphere.     
Suppose that the singular point set of $P$ consists of double points.
Then does it follow that $P$ is the projection of an unknotted $n$-knot? 
\enddefinition

In this paper we prove that the answer to Problem 1.2 in the case $n\geqq3$ is negative 
and hence the answer to Problem 1.1 in the case $n\geqq3$ is also negative.

\proclaim 
{Theorem 1.3}
Let $n$ be any integer greater than two.  
Then there exists a projection $P$ having the following properties.   
\roster
\item
$P$ is not the projection of any unknotted knot.   
\item
The singular point set of $P$ consists of double points.   
\item
$P$ is the projection of an $n$-knot which is diffeomorphic to the standard sphere.   
\endroster
\endproclaim 

In \S4 we prove 
there exists an immersed $n$-sphere ( $\subset\Bbb R^{n+1}\x\{0\}$ ) 
which is not the projection of any $n$-knot ($n>2$) (Theorem 4.1). 
Note that this theorem is different from Theorem 1.3.

\head \S 2. Proof of Theorem 1.3 in the case $n=3$ \endhead

Let $V$ be the punctured Kummer surface.   
It is known that $V$ has the following properties.  
 See, e.g., p.9 of 
[15]. 

\roster
\item
The intersection matrix is 
$\oplus^3$
$\pmatrix
0&1\\
1&0\\
\endpmatrix$
$\oplus^2$$E_8$.   
The signature $\sigma(V)$ is equal to $-16$.  
\item
There exists a handle decomposition with one 0-handle and twenty-two 2-handles.   
\endroster

We prove:
\proclaim{Lemma 2.1} 
There exist an immersion $\alpha:V\i$ $\Bbb R^4$ 
and charts $(U_i, \phi_i)$  of  $\Bbb R^4$ 
with the following properties. ($i=1,...,\mu$, where $\mu$ is a positive integer.  )

\hskip1cm
(i)
$U_1$ $\cup...\cup$ $U_\mu$ covers the singular point set of $\alpha(V)$.   

\hskip1cm
(ii)
$U_i\cap U_j=\phi$ for $i\neq j$.   

\hskip1cm
(iii)
$\phi_i$:$U_i\cong$ $\{(x,y,z,w)\vert$ $x^2+y^2+z^2+w^2<1\}$.   

\hskip1cm
(iv)
$\alpha^{-1}(U_i)$ consists of two components 
 $V_{1i}$ and $V_{2i}$ and 
$\phi_i\circ\alpha\vert_{V_{1i}}$ 
and 
$\phi_i\circ\alpha\vert_{V_{2i}}$ 

\hskip1.4cm
are embeddings.   
Furthermore:

\hskip17mm
$\phi_i\circ\alpha(V_{1i})$= $\{(x,y,z,w)\vert$ $x^2+y^2<1, 
z^2+w^2\leqq\frac{1}{4}\}$. 
\hskip17mm
$\phi_i\circ\alpha(V_{2i})$= $\{(x,y,z,w)\vert$ 
$x^2+y^2\leqq\frac{1}{4}, z^2+w^2<1\}$.   
\endproclaim

\f{\bf Proof. } 
Take the handle decomposition of $V$$=h^0\cup (\cup_{i=1}^{22}h^2_i)$
associated with the framed link given on p.9 of 
[15].  
Let $L_V=$$(K_{V,1},.  .  .  ,K_{V,22})$ denote this framed link.   
As described in 
[15],    
the framing of $K_{V,1}$ is zero 
and the framing of $K_{V,i}$ is $-2$ for $i\neq1$.   
The underlying link is also denoted by $L_V$.   
The framing of attaching $h^2_i$ is that of $K_{V,1}$.

Consider $\Bbb R^4$ as $\Bbb R^3\x\{t\vert t\in\Bbb R\}$.   
Take an embedded closed 4-ball $B^4$ in 
$\Bbb R^3\x\{t\vert t\in\Bbb R\}$  
so that 
$B^4\subset$ $\Bbb R^3\x\{t\vert t\leqq0\}$  
and 
$B^4\cap$ $(\Bbb R^3\x\{t\vert t=0\})$ is the 3-ball $D^3$.  
Take $L_V$ in Int $D^3$.

Recall the following. 
Suppose $\Bbb R^2_1$ and  $\Bbb R^2_2$ intersect transversely in $\Bbb R^4$. 
Suppose $\Bbb R^2_1$, $\Bbb R^2_2$, $\Bbb R^4$ are oriented. 
Suppose $\Bbb R^2_1\cap\Bbb R^2_2$ is one point $p$. 
Then we can define the signature $\sigma(p)$ for $p$ 
by using the orientations of $\Bbb R^2_1$, $\Bbb R^2_2$, $\Bbb R^4$.

\proclaim{Sublemma}
There exists a self-transverse immersion 
$f:D^2_1\amalg.  .  .  \amalg D^2_{22}\i$ $\Bbb R^3\x[0,\infty)$ 
with the following properties.   
\roster
\item
$f$ is transverse to $\Bbb R^3\x\{0\}$.   
$f(D^2_1\amalg...\amalg D^2_{22})\cap$ $\Bbb R^3\x\{0\}=$ 
$f(\partial D^2_1\amalg...\amalg \partial D^2_{22}$).    
\item 
$f\vert_{\partial D^2_1\amalg...\amalg \partial D^2_{22}}$ 
is an embedding.   
\item  
$f$($\partial D^2_1\amalg...\amalg \partial D^2_{22}$) 
in  $\Bbb R^3\x\{0\}$ is $L_V$.  
\item  
Let $p_{ij}$ be the singular points  of $f(D^2_i)$. 
We have  

$\Sigma_j\sigma(p_{ij})$
=$\cases
\text{0}  & \text{ for $i=1$} \\  
\text{$-1$} & \text{ for $i=2,.  .  ,22.  $} \\ 
\endcases$ 
\endroster
\endproclaim

Proof. Since $\pi_1(\Bbb R^3\x[0,1))=1$, 
there is a map $g:D^2_1\amalg.  .  .  \amalg D^2_{22}\i$ $\Bbb R^3\x[0,\infty)$ 
such that 
$g(\partial D^2_i)=K_{V,i}$. 
Perturb $g$ and make $g$ self-transverse. 
Let $q_{ik}$ be the singular points of $g(D^2_i)$.

We prove: 
There is a self-transverse immersed 2-sphere $S_\varepsilon\subset\Bbb R^3\x(0,\infty)$ 
 ($\varepsilon=+1,-1$) 
such that the singular point set is one point and its signature is $\varepsilon$.

Because: Take  $X=\{(a,b,c,d)\vert a^2+b^2+c^2+d^2\leqq1\}\subset\Bbb R^3\x(0,\infty)$. 
\nl Put $D^2_1=\{(a,b,c,d)\vert a^2+b^2\leqq1, c=d=0\}$.
Put $D^2_2=\{(a,b,c,d)\vert a=b=0, c^2+d^2\leqq1\}$.     
Then $H=(\partial D^2_1, \partial D^2_2)$ in the 3-sphere $\partial X$ is 
the Hopf link. 
Take a Seifert surface $A$ for $H$ which is diffeomorphic to $S^1\x[-1,1]$. 
Then $D_1\cup A\cup D_2$ is an example of $S^2_{+1}$. 
\nl Put $D^2_1=\{(a,b,c,d)\vert a^2+c^2\leqq1, b=d=0\}$ and 
$D^2_2=\{(a,b,c,d)\vert a=c=0, b^2+d^2\leqq1\}$ instead. 
Then we obtain $S^2_{-1}$. 
This completes the proof.

We continue the proof of Sublemma. 
Take a connected sum of $S^2_{\varepsilon}$ and $g(D^2_\mu)$. 
Then we obtain an immersed 2-disc $D$. 
Put $g':D^2_1\amalg...\amalg D^2_{22}\i$ $\Bbb R^3\x[0,\infty)$   
so that $g'(D^2_i)=g(D^2_i)$ for $i\neq\mu$ and  $g'(D^2_\mu)=D$. 
Let $q'_{ik'}$ be the singular points of $g'(D^2_i)$.
Then $\Sigma_{k'}\sigma(q'_{ik'})=\Sigma_k\sigma(q_{ik})$ for $i\neq\mu$ 
and $\Sigma_{k'}\sigma(q'_{\mu k'})=\Sigma_k\sigma(q_{\mu k})+\varepsilon$.

We make $f$ in Sublemma from the above $g$.  
Suppose $\Sigma_k\sigma(q_{ik})\neq\Sigma_j\sigma(p_{ij})$. 
By using immersed spheres as $S^2_{\varepsilon}$ and using operations as above, 
we modify $g$ to obtain $f$.  
This completes the proof of Sublemma.

Let $\nu_i$ be the normal $D^2$-bundle of $D^2_i$ in $\Bbb R^3\x[0,\infty)$. 
Let $E_i$ be the total space of $\nu_i$. 
Note $E_i\cong D^2_i\x D^2$. 
Here, we have an immersion $\beta_i:E_i\i\Bbb R^3\x\{0\}$. 
Note $\beta_i(E_i)\cap B^4=$ $(\partial D^2_i)\x D^2\subset \Bbb R^3\x\{0\}$. 
Let $W=B^4\cup(\cup_{i=1}^{22}E_i)$. 
Here, we have an immersion $\beta:W=B^4\cup(\cup_{i=1}^{22}E_i)\i\Bbb R^4$   
such that $\beta(B^4)=B^4$ and $\beta(E_i)=\beta_i(E_i)$. 
Note we can regard $\cup_{i=1}^{22}E_i$ as 4-dimensional 2-handles attached to $B^4$ 
along $L_V$. 
$E_i$ is called $h_i$. 
We write $W=h^0\cup(\cup_{i=1}^{22}h^2_i)$. 

We prove:
The framing of attaching $h^2_i$ is $2\cdot\Sigma_j\sigma(p_{ij})$. 

Proof. 
We perturb $D^2_i=D^2_i\x\{0\}$ in $h^2_i=E_i=D^2_i\x D^2$ 
so that we perturb $\partial D^2_i$ in $\partial h^2_i$. 
The result is called $\bar D^2_i$. 
Suppose $D^2_i$ and $\bar D^2_i$ intersect transversely. 
Suppose $(\partial D^2_i, \partial\bar D^2_i)$ in $\Bbb R^3\x\{0\}$ is the trivial link. 
By the definition of the framing, 
the algebraic intersection number of  $D^2_i\cap\bar D^2_i$ is the framing. 
Note that the points  $D^2_i\cap\bar D^2_i$ exist arround the points $ p_{ij}$. 
The algebraic intersection number is  $2\cdot\Sigma_j\sigma(p_{ij})$. 
This completes the proof.

Hence  $W=h^0\cup(\cup_{i=1}^{22}h^2_i)$ is $V$. 
This completes the proof of Lemma 2.1.  

We use  Lemma 2.1 and prove the following Proposition 2.2. 

By the definition of 3-knots, 
all 3-knots are diffeomorphic to the standard sphere.  
Hence Proposition 2.2 induces Theorem 1.3 in the case $n=3$.

\proclaim {Proposition 2.2} 
Let $r$ be any integer.   
Then there exists a projection $P$ having the following properties.   
\roster
\item
If the projection of $K$ is $P$, $K$ is not unknotted.
\item
If the projection of $K$ is $P$, $\sigma(K)=16r$.   
\item
The singular point set of $P$ consists of double points.  
\endroster
\endproclaim

{\bf Note.}  
An integer $\sigma$ is the signature of a 3-knot 
if and only if 
$\sigma$ is a multiple of sixteen.    
See, e.g., \S10 of 
[19].  

{\bf Proof of Proposition 2.2.}  We first prove Proposition 2.2 when $r=-1$.  

Consider 
$\alpha:V\i$ $\Bbb R^4$ in Lemma 2.1 as 
$\alpha:V\i$ $\Bbb R^4\x\{0\}$ 
($\subset\Bbb R^4\x\Bbb R$).   
Let $\beta$ be the self-transverse immersion 
$\alpha\vert_{\partial V=S^3}$:$S^3\i\Bbb R^4\x\{0\}$.   
We prove 

\proclaim{Lemma}
$\beta(S^3)$ is the projection of a 3-knot $(\subset\Bbb R^4\x\Bbb R)$. 
\endproclaim

Proof. 
Let $\pi:$
$\Bbb R^4\x\Bbb R$ $\rightarrow$ $\Bbb R^4\x\{0\}$
be the natural projection map.    
There is a submanifold $B$ which is diffeomorphic to $V$ 
such that $\pi(B)=\alpha(V)$.
Because: The immersion $\alpha$ is an embedding in 
$\{\Bbb R^4-(U_1\cup...\cup U_\mu)\}\x\Bbb R$. 
In each $U_i$, 
push $\alpha(V_{1i})$ into the direction of 
$\Bbb R^4\x\{t\vert t>0\}$ (or $\Bbb R^4\x\{t\vert t<0\}$ ) 
so that $\alpha(V_{2i})$ is fixed. 
Thus we obtain a submanifold $B$ which is diffeomorphic to $V$ 
such that $\pi(B)=\alpha(V)$.

Hence $\beta(S^3)$ is the projection of 
the 3-knot $\partial B$. 
This completes the proof of the Lemma.

Let $P$ denote the projection $\beta(S^3)$.    

\proclaim{Claim 2.3} 
$P$ satisfies conditions (1)(2)(3) of Proposition 2.2 when $r=-1$. 
\endproclaim 
  
{\bf Proof of Claim 2.3. }
Let $Q$ be the singular point set of $P$. 
By the construction of $\alpha$ and $\beta$, 
we have 

(1)
$Q\subset (U_1\cup...\cup U_\mu)\x\{t=0\}$ 

(2)$Q\cap U_i$ is $S^1\x S^1$. 
(In particular, see Lemma 2.1(iv). )
Hence $Q\cap U_i$ is connected. 

(3)$Q\cap U_i$ consists of double points.  

Hence $P$ satisfies condition (3) of Proposition 2.2.

Let $A$ be a 3-knot whose projection is $P$. 
Put $W_{1i}$=$(\pi\vert A)^{-1}(\partial V_{1i})$  
\nl and $W_{2i}$=$(\pi\vert A)^{-1}(\partial V_{2i})$. 
Recall that $P\supset \partial V_{1i}$ 
and $P\supset \partial V_{2i}$.  
Here, we write $V_{ji}$ for $\alpha(V_{ji})$.

We assign to the 3-knot $A$ an element $\rho(A)\in\Bbb Z_2^\mu$ given as follows. 
Let $\Bbb Z_2=\{+1, -1\}$. 
If $W_{1i}$ is over (resp. under ) $W_{2i}$, 
then we define the $i$-th coordinate of 
$\rho(A)$ to be $+1(resp. -1)$. 
We define  `over' and `under' by using the direction of $\Bbb R_t$.

Note the following.  
(1)Let $A'$ be a 3-knot whose projection is $P$. 
If $\rho(A')=\rho(A)$, then $A'$ is equivalent to $A$. 
(2)For any element $x\in\Bbb Z_2^\mu$, 
there is such a 3-knot $K$ with   $\rho(K)=x$. 

Let $B$ be a submanifold which is diffeomorphic to $V$ 
such that $\pi(B)=\alpha(V)$.
Put $X_{1i}$=$(\pi\vert A)^{-1}( V_{1i})$  
and $X_{2i}$=$(\pi\vert A)^{-1}( V_{2i})$. 
We give the submanifold $B$ 
an element $\rho(B)\in\Bbb Z_2^\mu$ as follows. 
If $X_{1i}$ is over (resp. under ) $X_{2i}$, 
then we define the $i$-th coordinate of $\rho(B)$ to be $+1(resp. -1)$. 

Note the following.  
(1)Let $B'$ be 
a submanifold which is diffeomorphic to $V$ 
such that $\pi(B')=\alpha(V)$.    
If   $\rho(B')=\rho(B)$, then 
the submanifold $B'$ is equivalent to the submanifold $B$. 
(2)For any element $x\in\Bbb Z_2^\mu$, 
there is such a submanifold $B$ with $\rho(B)=x$.

Take a 3-knot $A$ whose projection is $P$. 
Then there is a submanifold $B$ 
which is diffeomorphic to $V$ 
such that $\rho(B)=\rho(A)$. 
Since $\rho(B)=\rho(\partial B)$, $\rho(\partial B)=\rho(A)$. 
Hence the 3-knot $A$ is equivalent to the 3-knot $\partial B$. 
Therefore $A$ has a Seifert hypersurface 
which is diffeomorphic to $V$. 
Hence the signature of $A$ is $-16$.   
Therefore $A$ is knotted.  
Hence $P$ satisfies conditions (1) (2) of Proposition 2.2. 

This completes the proof of Claim 2.3 and 
thus the proof of Proposition 2.2 in the case $r=-1$. 
We next prove Proposition 2.2 when $r\neq -1$. 
We divided the proof into the three cases, $r< -1$, $r=0$, $r\geqq 1$.

Let $-P$ denote what we obtain from $P$ when we give the opposite orientation to $P$. 
Let $P^*$ denote what we obtain from $P$ 
when we give the opposite orientation to $S^3$.   


We prove Proposition 2.2 when $r<-1$. 
Let $\bar P$ be an immersed 3-sphere in  
$\Bbb R^4\x\{0\}=\Bbb R^3\x\Bbb R_w\x\{0\}$. 
Suppose 
$\bar P\subset$  $\Bbb R^3\x\{w\vert 0\leqq w\leqq\vert r\vert\}\x\{0\}$. 
For each $s\in\Lambda=\{0,1...,\vert r\vert\}$, suppose that 
$\bar P$ $\cap$ ($\Bbb R^3\x$ $\{w=s\}\x\{0\}$ ) 
is a 2-sphere which bounds a 3-disc $D_s$ embedded in 
$\Bbb R^3\x$ $\{ w=s\}$ $\x\{0\}$. 
If $s, s+1\in\Lambda$, put $P_{s+1}$=
$\bar P$ $\cap$ ($\Bbb R^3\x$ $\{s\leqq w\leqq s+1\}\x\{0\}$ ).  
Suppose $P_{s+1}\cup D^3_s\cup D^3_{s+1}$ is a parallel displacement of $P$. 
Then $\bar P$ is the projection of a 3-knot. 
So $\bar P$ satisfies condition (3) of Proposition 2.2. 
It also follows that the 3-knots with projection $\bar P$ are of the form 
$K_1$ $\sharp$... $\sharp$ $K_{\vert r\vert}$, 
where the projection of $K_*$ is $P$. 
Note that 
$\sigma$ ($K_1\sharp$... $\sharp K_{\vert r\vert}$)
= $\Sigma_{*=1}^{\vert r\vert}$ $\sigma(K_*)$ 
=$-16\vert r\vert$  
=16$r$. 
Hence $\bar P$ satisfies condition (2) of Proposition 2.2. 
Also each $K_1\sharp$... $\sharp K_{\vert r\vert}$ is not unknotted. 
Hence $\bar P$ satisfies condition (1) of Proposition 2.2. 
This completes the proof of Proposition 2.2 in the case of $r<-1$.

We now prove Proposition 2.2 when $r\geqq1$. 
Let $\bar P$ be as above. 
Take $-\bar P$.
Then $-\bar P$ satisfies the conditions of Proposition 2.2 
when $r\geqq1$. 

Finally, we prove Proposition 2.2 when $r=0$. 
Let $\hat P$ be an immersed 3-sphere in  $\Bbb R^4\x\{0\}=$ $\Bbb R^3\x\Bbb R_w\x\{0\}$.
Suppose $\hat P$ $\cap$ ($\Bbb R^3\x$ $\{w=0\}\x\{0\}$ ) 
is a 2-sphere which bounds a 3-disc $D$ embedded in 
$\Bbb R^3\x$ $\{w=0\}$ $\x\{0\}$. 
Suppose 
[$\hat P$ $\cap$ 
($\Bbb R^3\x$ $\{ w\geqq 0\}$$\x\{0\}$] 
$\cup D$ is a parallel displacement of $P$ 
and 
[$\hat P$ $\cap$ 
($\Bbb R^3\x$ $\{ w\leqq 0\}$$\x\{0\}$] 
$\cup D$ is a parallel displacement of $-P^*$.  
Then $\hat P$ is the projection of a 3-knot. 
So $\hat P$ satisfies condition (3) of Proposition 2.2. 
It also follows that all 3-knots with projection $\hat P$ 
are of the form $K_1\sharp (-K^*_2)$, 
where the projection of $K_*$ is $P$. 
Note that $\sigma$ ($K_1\sharp (-K^*_2)$) 
= $\sigma(K_1)+\sigma(-K^*_2)=$0.  
Hence $\hat P$ satisfies condition (2) of Proposition 2.2.

In order to prove that $\hat P$ satisfies condition (1) of Proposition 2.2, 
we prove that each $K_1\sharp (-K^*_2)$ is not unknotted.  
We begin by recalling the following fact. 
See, e.g., \S 14 of 
[19]  
 for the Alexander polynomials.   
See, e.g., \S 6 of 
[19] 
 for simple knots.   
The author gives a proof in the appendix.

\proclaim {Theorem 2.4 (known)}    
Let $K$ be a simple $(2k+1)$-knot $(k\geq1)$.   
Then the Alexander polynomial of $K$ is trivial if and only if $K$ is unknotted.   
\endproclaim

 $K_1$ and $K_2$ are 3-knots whose projections are $P$. 
Then $K_i$ bounds a Seifert hypersurface diffeomorphic to 
the punctured Kummer surface $V$. 
Hence $K_i$ bounds a simply connected Seifert hypersurface.  
So $K_i$ is a simple knot.  
Since $\sigma(V)\neq0$, $K_1$ is not unknotted and $-K_2^*$ is not unknotted.  
Therefore the Alexander polynomial 
$\Delta_{K_1}(t)$ of $K_2$ is nontrivial and the Alexander polynomial 
$\Delta_{-K_2^*}(t)$ of $-K_2^*$ is also nontrivial.

Recall the following fact.   

\proclaim {Theorem (known)}    
Let $K_1$ and $K_2$ be $(2k+1)$-knots $(k\geq1)$.   
Suppose the Alexander polynomial of each $K_i$ is nontrivial.   
Then the Alexander polynomial of $K_1\sharp K_2$ is also nontrivial.   
\endproclaim

From this we conclude that each $K_1\sharp (-K^*_2)$ is not unknotted. 

This proves condition (1) of Proposition 2.2 in the case $r=0$ and 
thus completes the proof of Proposition 2.2. 
We conclude that Theorem 1.3 holds in the case of $n=3$.

\head \S 3. Proof of Theorem 1.3 in the case $n>3$ \endhead

We define the {\it spun projection} of a projection. 
Let $P$ be the projection of an $n$-knot $K$. 
Suppose $K\subset$ $\Bbb R^{n}\x\Bbb R_t\x\Bbb R_u$  
and $P\subset$ $\Bbb R^{n}\x  \{t=0\}\x\Bbb R_u$. 
Let $\pi:\Bbb R^{n}\x\Bbb R_t\x\Bbb R_u\rightarrow$  
$\Bbb R^{n}\x  \{t=0\}\x\Bbb R_u$ 
be the natural projection map. 
We suppose 
$P\subset$ $\Bbb R^{n}\x  \{t=0\}\x\{u\geqq0\}$ and that 
$P\cap(\Bbb R^{n}\x  \{t=0\}\x\{u=0\})$  is an $n$-disc $D$ embedded in 
$\Bbb R^{n}\x  \{t=0\}\x\{u=0\}$. 
Suppose $D$ does not intersect the singular point set $Q$ of $P$.

Take  $(\Bbb R^{n}\x  \Bbb R_t\x\Bbb R_u)\x\Bbb R_v$ and we regard 
$\Bbb R^{n}\x  \Bbb R_t\x\Bbb R_u$ 
as $\Bbb R^{n}\x  \Bbb R_t\x\Bbb R_u \x\{v=0\}$. 
We regard  
$\Bbb R^{n}\x  \Bbb R_t\x\Bbb R_u\x\Bbb R_v$ 
as the result of rotating 
$\Bbb R^{n}\x  \Bbb R_t\x \{u\geqq0\}\x\{v=0\}$ 
\nl around the axis 
$\Bbb R^{n}\x  \Bbb R_t\x \{u=0\}\x\{v=0\}$. 

Then 
we regard  
$\Bbb R^{n}\x  \{t=0\}\x\Bbb R_u\x\Bbb R_v$ 
as the result of rotating 
\nl $\Bbb R^{n}\x  \{t=0\}\x \{u\geqq0\}\x\{v=0\}$ 
around the axis 
$\Bbb R^{n}\x  \{t=0\}\x \{u=0\}\x\{v=0\}$. 

Let 
$\pi:\Bbb R^{n}\x\Bbb R_t\x\Bbb R_u\x\Bbb R_v\rightarrow$  
 $\Bbb R^{n}\x  \{t=0\}\x\Bbb R_u\x\Bbb R_v$ 
be the natural projection map. 
We can regard $P$ as a subset of $\Bbb R^{n}\x\{t=0\}\x\{u\geqq0\}\x\{v=0\}$ 
\nl since 
$\Bbb R^{n}\x\{t=0\}\x\{u\geqq0\}$  
is identified with  
$\Bbb R^{n}\x\{t=0\}\x\{u\geqq0\}\x\{v=0\}$.

When we rotate  
$\Bbb R^{n}\x  \{t=0\}\x \{u\geqq0\}\x\{v=0\}$ 
around the axis 
\nl $\Bbb R^{n}\x  \{t=0\}\x \{u=0\}\x\{v=0\}$,  
we rotate $\overline{P-D}$ as well. 
The result, denoted $\widetilde P$, 
 is called  the {\it spun projection} of $P$.

Let $\widetilde K$ be an $(n+1)$-knot.
Let $\widetilde K$ be the spun knot of an $n$-knot $K$. 
See 
[26]  
for a basic description of spun knots. 
Then we can suppose the following (1)(2). 

(1)$K$ is in 
$\Bbb R^{n}\x  \Bbb R_t\x \{u\geqq0\}\x\{v=0\}$. 

(2)$K\cap$  
($\Bbb R^{n}\x  \Bbb R_t\x \{u=0\}\x\{v=0\}$) 
is the $n$-disc $D$ which is defined above.

When we rotate 
$\Bbb R^{n}\x  \Bbb R_t\x \{u\geqq0\}\x\{v=0\}$ 
around the axis 
\nl$\Bbb R^{n}\x  \Bbb R_t\x \{u=0\}\x\{v=0\}$, 
 we rotate $\overline{K-D}$ as well. 
The result is denoted $\widetilde K$.

\proclaim {Lemma 3.1}
The above  $\widetilde P$ is the projection of $\widetilde K$. 
\endproclaim 

{\bf Proof. }  
Let 
$\alpha_\theta:\Bbb R^{n}\x  \Bbb R_t\x\Bbb R_u\x\Bbb R_v$ 
$\rightarrow\Bbb R^{n}\x  \Bbb R_t\x\Bbb R_u\x\Bbb R_v$ 
be the map 
\nl $(*,*,u,v)\mapsto$ 
$(*,*,u\cdot cos\theta-v\cdot sin\theta, $
$u\cdot sin\theta+v\cdot cos\theta)$. 
\nl We can regard $\pi$ as 
$\widetilde\pi\vert$
$(\Bbb R^{n}\x\Bbb R_t\x\{u=0\}\x\{v=0\})$. 
\nl Then 
$\widetilde\pi\vert$
$(\Bbb R^{n}\x  
\Bbb R_t\x\{u=r\cdot cos\theta, 
v=r\cdot sin\theta\vert r\geqq0, 
\theta$ is fixed$\})$
is equal to the compositions 
($\alpha_\theta)\circ\widetilde\pi\circ$($\alpha_\theta^{-1})$. 
\nl Furthermore 
$\widetilde K\cap$
$(\Bbb R^{n}\x  
\Bbb R_t\x\{u=r\cdot cos\theta, 
v=r\cdot sin\theta\vert r\geqq0, 
\theta$ is fixed$\})$
=$\alpha_\theta(\overline{K-D})$ 
\nl and 
$\widetilde P\cap$
$(\Bbb R^{n}\x  
\Bbb R_t\x\{u=r\cdot cos\theta, 
v=r\cdot sin\theta\vert r\geqq0, 
\theta$ 
is fixed$\})$
=$\alpha_\theta(\overline{P-D})$. 
\nl Hence $\widetilde\pi(\widetilde K)=\widetilde P$. 
This completes the proof.

Let $Q$ be the singular point set of $P$. 
Let $\widetilde Q$ be the singular point set of $\widetilde{P}$. 
Then the following holds.

\proclaim{Lemma 3.2} 
When we rotate  
$\Bbb R^{n}\x  \{t=0\}\x \{u\geqq0\}\x\{v=0\}$ 
around the axis 
\nl$\Bbb R^{n}\x  \{t=0\}\x \{u=0\}\x\{v=0\}$ 
to make $\widetilde P$ from $P$, 
we rotate $Q$ as well.  
Then the result $\widetilde Q$ coincides with $Q\x S^1$ and 
$\widetilde Q$ consists of double points. 
\endproclaim

We use the above preliminaries to prove Theorem 1.3 when $n>3$.


Let $P^{(3)}$ be the projection $\hat P$ of a 3-knot,  
which is defined in the proof of the $r=0$ case of Proposition 2.2.  
Let $P^{(n+1)}$ be the spun projection of $P^{(n)}$ ($n\geqq3$).

\proclaim{Claim 3.3}
$P^{(n)}$ satisfies conditions (1)(2)(3) of Theorem 1.3. ($n\geqq3$). 
\endproclaim

\f{\bf Proof.} 
By Proposition 2.2, there is a 3-knot $K^{(3)}$ 
which is diffeomorphic to the standard sphere 
and whose projection is $P^{(3)}$.  
Suppose $K^{(n+1)}$ is the spun knot of $K^{(n)}$ ($n\geqq3$). 
By induction on $n$, $K^{(n)}$ is diffeomorphic to the standard sphere($n\geqq3$).  
By Lemma 3.1, $P^{(n)}$ is the projection of  $K^{(n)}$. 
Hence $P^{(n)}$ satisfies condition (3) of Theorem 1.3.

Let $Q^{(n)}$ be the singular point set of $P^{(n)}$.  
By Lemma 3.2,  
$Q^{(n)}$ consists of double points. 
Hence $P^{(n)}$ satisfies condition (2) of Theorem 1.3.

In order to prove $P^{(n)}$ satisfies condition (1) of Theorem 1.3, 
we prove: 
\proclaim{Lemma 3.4} 
Let $X$ be an $(n+1)$-knot whose projection is $P^{(n+1)}$. 
Then there is an $n$-knot $Y$ such that 
the projection of $Y$ is  $P^{(n)}$ and that 
$X$ is the spun knot of $Y$. 
 ($n\geqq3$.)
\endproclaim

{\bf Proof. }
By Lemma 3.2, 
The number of the connected components of $Q^{(n)}$ is $\mu$. 
Let $Q_1^{(n)}, ...,Q_\mu^{(n)}$ be the connected components of $Q^{(n)}$. 
Suppose $Q_i^{(n+1)}$ is made from $Q_i^{(n)}$ ( $n\geqq3$ ). 

Let $A$ be an $n$-knot whose projection is $P$. 
Let $Q_{1i}^{(n)}$, $Q_{2i}^{(n)}$ be the connected components of 
$(\pi\vert A)^{-1}(Q_i^{(n)})$. 
Note that 
$Q_{1i}^{(n)}\subset A$ and  $Q_{2i}^{(n)}\subset A$. 
We assign to the $n$-knot $A$ 
an element $\rho(A)\in\Bbb Z_2^\mu$ as follows. 
Let $\Bbb Z_2=\{+1, -1\}$. 
If $Q_{1i}$ is over (resp. under ) $Q_{2i}$, 
then we define the $i$-th coordinate of 
$\rho(A)$ to be $+1(resp. -1)$. 
We define  `over' and `under' by using 
the direction of $\Bbb R_t$. 
Note the following 
(1)Let $A'$ be another $n$-knot whose projection is $P$. 
If $\rho(A')=\rho(A)$, then $A'$ is equivalent to $A$. 
(2)For any element $x\in\Bbb Z_2^\mu$, 
there is such an $n$-knot $K$ 
with   $\rho(K)=x$.

Let $X$ be an $(n+1)$-knot whose projection is $P^{(n+1)}$. 
Then there is an $n$-knot $Y$ whose projection is   $P^{(n)}$ 
such that $\rho(Y)=\rho(X)$. 
Let $\widetilde Y$ be an $(n+1)$-knot which is
 the spun knot of the $n$-knot $X$. 
By Lemma 3.1, the projection of $\widetilde Y$ is $P^{(n+1)}$. 
Then $\rho(\widetilde Y)=\rho(Y)$.
Hence $\rho(\widetilde Y)=\rho(X)$, so the $(n+1)$-knot $X$ 
is equivalent to the $(n+1)$-knot $\widetilde Y$. 
By induction on $n$, Lemma 3.4 holds.

We continue the proof of Claim 3.3. 
If Lemma 3.5 below is true, 
then 
$P^{(n)}$ satisfies condition (1) of Theorem 1.3.

\proclaim{Lemma 3.5} 
Let $K$ be an $n$-knot whose projection is $P^{(n)}$. 
Then $K$ is not unknotted. 
$(n\geqq3.)$    
\endproclaim

\f{\bf Proof.  }
The $n=3$ case holds by \S2. We prove $n>3$. 
By Lemma 3.4, there are 
a 3-knot $Z^{(3)}$,   
a 4-knot $Z^{(4)}$,..., 
and an $n$-knot $Z^{(n)}$ such that 
$Z^{(q+1)}$ is the spun knot of $Z^{(q)}$ ($q=3,...,n-1$),  
that $Z^{(n)}=K$, 
and that the projection of $Z^{(r)}$ is $P^{(r)}$ ($r=3,...,n$).

Recall that the following facts hold by Theorem 4,1 of 
[11]  
or by using the Mayer-Vietoris exact sequence.    
See, e.g., \S 14 of 
[19]  
 for the Alexander polynomials.   
See, e.g., p.160 of 
[21] 
 and 
[18]  
for the Alexander invariant.   
Let $\widetilde {X}_{K}$ denote 
the canonical infinite cyclic covering 
of the complement of the knot $K$.

\proclaim {Theorem 3.6 (known)}    
Let $K$ be a simple $(2k+1)$-knot $(k\geq1)$.   
Let $\Delta_K(t)$ be the Alexander polynomial of $K$.   
Suppose  the $(k+1)$- Alexander invariant  
$H_{k+1}(\widetilde {X}_K;\Bbb Q)\cong$ 
 $\{\Bbb Q[t, t^{-1}]$/ $\delta^1_K(t)\}$  
 $\oplus .  .  .   \oplus$
 $\{\Bbb Q[t, t^{-1}]$/ $\delta^p_K(t)\}$.   
 Then  
 $\Delta_K(t)=$  
 $a\cdot t^b\cdot\delta^1_K(t)$ $\cdot.  .  .  \cdot$ $\delta^p_K(t)$ 
  for a rational number $a$ and an integer   $b$  
and we can put $\Delta_K(1)=1.$
\endproclaim

\proclaim {Theorem 3.7 (known) }   
Let $K^{(n+1)}$ be the spun knot of $K^{(n)}$ ($n\geqq1$).   
Let $H_{k}(\widetilde {X}_{K^{(n)}};\Bbb Q)$ 
(resp. $H_{k}(\widetilde {X}_{K^{(n+1)}};\Bbb Q)$  )
denote the $k$-Alexander invariant of $K^{(n)}$ (resp. $K^{(n+1)}$ ).   
Suppose that 
$K^{(n)}$ bounds a Seifert hypersurface $V$ such that 
$H_1(V;\Bbb Z)\cong0$.
Then 
$H_2(\widetilde {X^{(n+1)}};\Bbb Q)\cong$ 
$H_2(\widetilde {X^{(n)}};\Bbb Q)$. 
\endproclaim

\proclaim {Proposition 3.8 (known) }   
Let $K^{(n+1)}$ be the spun knot of $K^{(n)}$ ($n\geqq1$).   
If $K^{(n)}$ has a  simply connected Seifert hypersurface, 
then $K^{(n+1)}$ has a  simply connected 
Seifert hypersurface.  
\endproclaim

As in the proof of the $r=0$ case of Proposition 2.2, 
the Alexander polynomial $\Delta_{Z^{(3)}}(t)$ of $Z^{(3)}$ is nontrivial. 
By Theorem 3.6, the 2-Alexander invariant     
$H_2(\widetilde {X}_{Z^{(3)}};\Bbb Q)$  is nontrivial.                             
By Theorem 3.7 and Proposition 3.8, the 2-Alexander invariant     
$H_2(\widetilde {X}_{Z^{(s)}};\Bbb Q)$ is nontrivial ($s=4,..n$).                     

Therefore this also completes the proof of Lemma 3.5. 
This completes the proof Theorem 1.3 in the final case when $n>3$.

\head  \S 4. The proof of Theorem 4.1 
\endhead

An {\it immersed $n$-sphere $A$} is an image of a self-transverse immersion $S^n\i\Bbb R^{n+1}$.  
Let 
$\pi:\Bbb R^{n+1}\x\Bbb R$ $\to$ $\Bbb R^{n+1}\x\Bbb \{0\}$ 
be the natural projection map.    
We regard an immersed $n$-sphere $A$ as in $\Bbb R^{n+1}\x\{0\}$.  
We say that $A$ {\it lift}s (into $\Bbb R^{n+1}\x\Bbb R$ )  
if there is an $n$-knot $K$ whose projection is $A$. 
Then $K$ is called a {\it lift}. 

It is natural to ask the following question.   

{\bf Problem.}
Do all immersed $n$-spheres  ($\subset \Bbb R^{n+1}\x\{0\}$) lift into  $\Bbb R^{n+1}\x\Bbb R$?

It is well-known that, when $n=1$, the answer is affirmative. 
[10] 
proved: when $n=2$, the answer is negative. 
( [6] 
gave an alternative proof.)  

In this section, we give the negative answer to the above Problem in the case $n>2$.

\proclaim{ Theorem 4.1 } 
Let $n$ be any integer greater than two. 
There exists an immersed $n$-sphere ($\subset\Bbb R^{n+1}\x\{0\}$)
which does not lift into $\Bbb R^{n+1}\x\Bbb R$.  
\endproclaim

{\bf Note.} Theorem 4.1 is different from Theorem 1.3.

\f{\bf Proof of Theorem 4.1.}
We define the {\it spun immersed $(n+1)$-sphere} of an immersed $n$-sphere. 
Let $A$ be an immersed $n$-sphere. 
Take $\Bbb R^{n}\x\{x\vert x\in\Bbb R\}\x\{y\vert y\in\Bbb R\}$.  
We suppose $A\subset \Bbb R^{n}\x\{x\geqq0\}\x\{y=0\}$.   
We suppose $A\cap(\Bbb R^{n}\x  \{x=0\}\x\{y=0\})$  
is an $n$-disc $D$ embedded in $\Bbb R^{n}\x  \{x=0\}\x\{y=0\}$. 
We suppose $D$ does not intersect the singular point set of $A$. 
We regard  
$\Bbb R^{n}\x\{x\vert x\in\Bbb R\}\x\{y\vert y\in\Bbb R\}$ 
as the result of rotating 
\nl $\Bbb R^{n}\x\{x\geqq0\}\x\{y=0\}$ 
 around the axis  $\Bbb R^{n}\x  \Bbb R_t\x \{x=0\}\x\{y=0\}$. 
 When we rotate  
\nl$\Bbb R^{n}\x\{x\geqq0\}\x\{y=0\}$ 
around the axis $\Bbb R^{n}\x  \Bbb R_t\x \{x=0\}\x\{y=0\}$, 
 we rotate $\overline{A-D}$ as well. 
The result, denoted $\widetilde A$, 
 is called  the {\it spun immersed $(n+1)$-sphere} of $A$. 

We regard 
$\Bbb R^{n}\x\{x\vert x\in\Bbb R\}\x\{y\vert y\in\Bbb R\}$ 
as 

\f$\Bbb R^{n}\x\{x\vert x\in\Bbb R\}\x\{y\vert y\in\Bbb R\}\x\{z\vert z=0\}$ 
$\subset\Bbb R^{n}\x\{x\vert x\in\Bbb R\}\x\{y\vert y\in\Bbb R\}\x\{z\vert z\in\Bbb R\}$. 

We consider whether 
$\widetilde A$ ($\subset\Bbb R^{n}\x\{x\vert x\in\Bbb R\}\x\{y\vert y\in\Bbb R\}\x\{z=0\}$ )

\f lifts into 
$\Bbb R^{n}\x\{x\vert x\in\Bbb R\}\x\{y\vert y\in\Bbb R\}\x\{z\vert z\in\Bbb R\}$. 

We consider whether 
$A$ ($\subset\Bbb R^{n}\x\{x\geqq0\}\x\{y=0\}\x\{z=0\}$)  

\f lifts into 
$\Bbb R^{n}\x\{x\geqq0\}\x\{y=0\}\x\{z\vert z\in\Bbb R\}$.

We prove Claim 1. By the above result of 
[10]  
and Claim 1, Theorem 4.1 holds. 

\proclaim{Claim 1} 
Let $\widetilde A$ be the spun immersed $(n+1)$-sphere of an immersed $n$-sphere $A$. 
If $\widetilde A$ lifts, then $A$ lifts. 
\endproclaim

\f{\bf Proof of Claim 1. }
Let $E$ be an immersed $m$-sphere in $\Bbb R^{m+1}\x\{0\}$. 
Suppose $E$ lifts into $\Bbb R^{m+2}=\Bbb R^{m+1}\x\Bbb R$. 
Let $K$ be a lift of $E$. 
Let $Q$ be a compact $n$-submanifold $\subset E$. 
Suppose $Q$ does not intersect with the singular point set of $E$. 
By using the partition of unity, we can suppose that $E\cap K=Q$. 

Let $\widetilde K$ 
$\subset \Bbb R^{n}\x\{x\vert x\in\Bbb R\}\x\{y\vert y\in\Bbb R\}\x\{z\vert z\in\Bbb R\}$ 

\f be a lift of  
 $\widetilde A$ 
$\subset \Bbb R^{n}\x\{x\vert x\in\Bbb R\}\x\{y\vert y\in\Bbb R\}\x\{z=0\}$. 
We suppose 
$\widetilde A\cap \widetilde K=\partial D$.   

Let $K_0$=
$\widetilde K \cap(\Bbb R^{n}\x\{x\geqq0\}\x\{y=0\}\x\{z\vert z\in\Bbb R\})$. 
Then $K_0\cap D$ is $\partial D=\partial K_0$.

\f Then $K=K_0\cup D$ is an $n$-submanifold 
$\subset\Bbb R^{n}\x\{x\geqq0\}\x\{y=0\}\x\{z\vert z\in\Bbb R\}$.  

\f Then the projection of $K$ is $A$. 
This completes the proof of Claim 1 and Theorem 4.1.

\proclaim{ Theorem 4.2} 
Let $n>2$. 
Let $M$ be a closed manifold 
which can be immersed into $\Bbb R^{n+1}\x\{0\}$ and 
which can be embedded in $\Bbb R^{n+1}\x\Bbb R$.  
Then there exists an immersion $f:M\i \Bbb R^{n+1}\x\{0\}$ such that 
$f(M)$ does not lift into $\Bbb R^{n+1}\x\Bbb R$.  
\endproclaim 

\f{\bf Proof.} Take an immersion $g:M\i \Bbb R^{n+1}\x\{0\}$. 
Let $A$ be an immersed $n$-sphere  
which does not lift into $\Bbb R^{n+1}\x\Bbb R$.  
Take a connected sum $g(M)\sharp A$ in $\Bbb R^{n+1}\x\{0\}$.
If $g(M)\sharp A$ lifts, then $A$ lifts. 
Hence  $g(M)\sharp A$ does not lift. 
Take $f$ so that $f(M)=g(M)\sharp A$.

\head Appendix.  
The proof of Theorem 2.4. 
\endhead

We review the definition of the Alexander polynomial 
$\Delta_K(t)$ for a $(2k+1)$-knot $K$.   
(See \S 14 of [19]. )        
Let $A$ be a Seifert matrix for $K$. 
We define  $\Delta_K(t)$ 
to be det$(A-(-1)^kA')$,  
where $A'$ is the transposed matrix.
Note that we identify 
$\Delta_K(t)$ 
with 
$(-1)^rt^s\Delta_K(t)$(for any $r,s\in\Bbb Z$).

\proclaim {Theorem 2.4 (known)}  
Let $K$ be a simple $(2k+1)$-knot $(k\geq1)$. 
Then the following conditions are equivalent. 
\roster
\item
The Alexander polynomial of the simple $(2k+1)$-knot $K$ is trivial. 
\item
The simple $(2k+1)$-knot $K$ is trivial.
\endroster
\endproclaim

\f{\bf Note.}
See the definition of simple knots  for \S 6 of [19].       

\f{\bf Proof. }
(2) $\Rightarrow$ (1) is obvious.
We prove (1) $\Rightarrow$ (2).

Let ${\widetilde X}$ be the infinite cyclic covering of the complement $X$ of 
the knot $K$. 
By (1) 
we have: $\pi_1(X)\cong\Bbb Z$, 
 $\pi_1({\widetilde X})\cong1$, and 
$H_i({\widetilde X})\cong0$  $(i>0)$.

By Hurewictz's theorem, we have:
$\pi_i({\widetilde X})\cong0$$(i>1).$

Therefore 
$\pi_i(X)\cong0$ $(i>1)$.

The following theorem is proved essentially in 
 [17],   
[22],  
[24]. 
By this theorem,  $K$ is the trivial knot. This completes the proof.

\proclaim {Theorem }([17][22][24])  
Let $K_1$ be an $n$-knot $(n\geq3)$ in $S^{n+2}$. 
Then $K$ is trivial if and only if the following conditions hold. 

(1)$\pi_1(S^{n+2}-K_1)\cong\Bbb Z$.   

(2)$\pi_i(S^{n+2}-K_1)\cong0$ $(i>1)$. 
\endproclaim

\enddocument